\documentclass[12pt]{amsart}
\usepackage{graphicx}
\usepackage{times}
\newtheorem{theorem}{Theorem}[section]
\newtheorem{proposition}[theorem]{Proposition}
\newtheorem{lemma}[theorem]{Lemma}

\theoremstyle{definition}
\newtheorem{definition}[theorem]{Definition}
\newtheorem{example}[theorem]{Example}

\newcommand{\Mg}{\ensuremath{\mathcal{M}_{g,1}}}
\newcommand{\Cg}{\ensuremath{\mathcal{C}_{g,1}}}
\newcommand{\CgI}{\ensuremath{\mathcal{C}_{g,1}^{\mathrm{irr}}}}
\newcommand{\Hg}{\ensuremath{\mathcal{H}_{g,1}}}

\newcommand{\Sg}{\ensuremath{\Sigma_{g,1}}}

\newcommand{\Hom}{\mathop{\mathrm{Hom}}\nolimits}

\newcommand{\Id}{\mathop{\mathrm{id}}\nolimits}
\newcommand{\Q}{\ensuremath{\mathbb{Q}}}
\newcommand{\Z}{\ensuremath{\mathbb{Z}}}

\newcommand{\PZ}{\ensuremath{\mathbb{Z}_{\ge 0}}}

\theoremstyle{remark}

\numberwithin{equation}{section}

\begin{document}

\title{Abelian quotients of monoids of homology cylinders} 

\date{\today}

\author{Hiroshi Goda}
\address{Department of Mathematics,
Tokyo University of Agriculture and Technology,
2-24-16 Naka-cho, Koganei,
Tokyo 184-8588, Japan}
\email{goda@cc.tuat.ac.jp}
\author{Takuya Sakasai}
\address{Department of Mathematics, 
Tokyo Institute of Technology, 
2-12-1 Oh-okayama, Meguro-ku, 
Tokyo 152-8551, Japan}
\email{sakasai@math.titech.ac.jp}

\thanks{
The authors are partially supported
by Grant-in-Aid for Scientific Research,
(No.~21540071 and No.~21740044),
Ministry of Education, Science, 
Sports and Technology, Japan. 
}

\keywords{Homology cylinder, homology cobordism, sutured manifold,  
sutured Floer homology}


\begin{abstract}
A homology cylinder over a surface 
consists of a homology cobordism between 
two copies of the surface and markings of its boundary. 
The set of isomorphism classes of homology cylinders 
over a fixed surface has a natural monoid structure and 
it is known that this monoid can be seen as an enlargement of 
the mapping class group of the surface. 
We now focus on abelian quotients of this monoid. 
We show that both the monoid of 
all homology cylinders and that of irreducible homology cylinders 
are not finitely generated and moreover they have big abelian quotients. 
These properties contrast with the fact that the mapping class group 
is perfect in general. 
The proof is given by applying sutured Floer homology theory 
to homologically fibered knots studied in a previous paper. 
\end{abstract}

\maketitle

\section{Introduction}\label{sec:intro}

Let $\Sg$ be a compact oriented connected surface of genus $g$ with 
one boundary component. {\it Homology cylinders} over $\Sg$, 
each of which consists of a homology cobordism $M$ 
between two copies of $\Sg$ and markings of 
both sides of the boundary of $M$, 
appeared in the context of the theory of finite type 
invariants for 3-manifolds (see Goussarov \cite{gou}, 
Habiro \cite{habiro}, Garoufalidis-Levine \cite{gl} and 
Levine \cite{levin}), and play an important role in a systematic 
study of the set of 3-manifolds. 
In our previous paper \cite{gs}, 
we observed their relationship to knot theory 
by introducing {\it homologically fibered knots}. 

The set $\Cg$ of isomorphism classes of homology cylinders 
over $\Sg$ becomes a monoid by the natural stacking operation. 
It is known that the monoid $\Cg$ contains 
the mapping class group $\Mg$ of 
$\Sg$ as the group of units (see Example \ref{ex:mgtocg}). 
Moreover, many techniques and invariants to study 
$\Mg$ such as Johnson homomorphisms 
and the Magnus representation 
can be extended to $\Cg$ (see \cite{gl}, \cite{sakasai}). 
By using them, we can observe that 
$\Cg$ and $\Mg$ hold many properties in common. 

Taking account of the similarity between $\Cg$ and $\Mg$, 
we now pay our attention to abelian quotients of them. 
It is known that $\Mg$ is a perfect group for $g \ge 3$, 
namely it has no non-trivial abelian quotients. 
In this paper, however, we will show 
that the opposite holds for $\Cg$. 

The outline of this paper is as follows. 
After introducing homology cylinders in Section 
\ref{section:cylinder}, we see that $\Cg$ has a big abelian quotient 
arising from the reducibility of a homology cylinder as a 3-manifold 
(Theorem \ref{prop:notFG1}). However this fact seems 
not to be suitable for our purpose of comparing $\Mg$ and $\Cg$ because 
all homology cylinders coming from $\Mg$ have irreducible underlying 
3-manifolds. Therefore we 
shall introduce the submonoid $\CgI$ of $\Cg$ consisting of 
{\it irreducible homology cylinders}. 
The main result is Theorem \ref{thm:notFG2} that $\CgI$ 
has a big abelian quotient originating in 
a quite different context 
from that of $\Cg$ mentioned above. 
In fact, we prove it in Section \ref{sec:proof} 
as an application of sutured Floer homology theory. 
Sutured Floer homology was defined first in \cite{juhasz}
by Juh\'asz, then an alternative definition was given 
in \cite{ni} by Ni.  
It is a variant of 
Heegaard Floer homology theory 
defined by Ozsv\'ath and Szab\'o 
(here we only refer to \cite{os}, which 
contains the results we use later, for details). 
In the last section, we discuss our results from the viewpoint of 
homology cobordisms of homology cylinders. 

The authors would like to thank Dr.~Motoo Tange for giving them 
a lecture about sutured Floer homology theory. They also would like to 
thank the referee for his/her 
careful reading of the manuscript and valuable comments. 
The final publication is available at www.springerlink.com.

\section{Homology cylinders}\label{section:cylinder}

We first recall the definition of homology cylinders, following 
Garoufalidis-Levine \cite{gl} and Levine \cite{levin}. 

\begin{definition}
A {\it homology cylinder\/}  $(M,i_{+},i_{-})$ {\it over} $\Sigma_{g,1}$ 
consists of a compact oriented 3-manifold $M$ 
with two embeddings $i_{+}, i_{-}: \Sigma_{g,1} \hookrightarrow \partial M$ 
such that:
\begin{enumerate}
\renewcommand{\labelenumi}{(\roman{enumi})}
\item
$i_{+}$ is orientation-preserving and $i_{-}$ is orientation-reversing; 
\item 
$\partial M=i_{+}(\Sigma_{g,1})\cup i_{-}(\Sigma_{g,1})$ and 
$i_{+}(\Sigma_{g,1})\cap i_{-}(\Sigma_{g,1})
=i_{+}(\partial\Sigma_{g,1})=i_{-}(\partial\Sigma_{g,1})$;
\item
$i_{+}|_{\partial \Sigma_{g,1}}=i_{-}|_{\partial \Sigma_{g,1}}$; and 
\item
$i_{+},i_{-} : H_{*}(\Sigma_{g,1};\mathbb Z)\to H_{*}(M;\mathbb Z)$ 
are isomorphisms. 
\end{enumerate}
\end{definition}

Two homology cylinders $(M,i_+,i_-)$ and $(N,j_+,j_-)$ over $\Sg$ 
are said to be {\it isomorphic} if there exists 
an orientation-preserving diffeomorphism $f:M \xrightarrow{\cong} N$ 
satisfying $j_+ = f \circ i_+$ and $j_- = f \circ i_-$. 
We denote by $\mathcal{C}_{g,1}$ the set of all isomorphism classes 
of homology cylinders over $\Sg$. We define a product operation on $\Cg$ by 
\[(M,i_+,i_-) \cdot (N,j_+,j_-)
:=(M \cup_{i_- \circ (j_+)^{-1}} N, i_+,j_-)\]
for $(M,i_+,i_-)$, $(N,j_+,j_-) \in \Cg$, 
so that $\mathcal{C}_{g,1}$ becomes a monoid with the 
identity element 
$(\Sigma_{g,1} \times [0,1], \mathrm{id} \times 1, \mathrm{id} \times 0)$. 

\begin{example}\label{ex:mgtocg}
For each diffeomorphism $\varphi$ of 
$\Sigma_{g,1}$ which fixes $\partial \Sigma_{g,1}$ pointwise, 
we can construct a homology cylinder by setting 
\[(\Sigma_{g,1} \times [0,1], \mathrm{id} \times 1, 
\varphi \times 0),\]
where collars of $i_+ (\Sigma_{g,1})$ and $i_- (\Sigma_{g,1})$ 
are stretched half-way along $(\partial \Sigma_{g,1}) \times [0,1]$. 
It is easily checked that the isomorphism class of 
$(\Sigma_{g,1} \times [0,1], \mathrm{id} \times 1, \varphi \times 0)$ 
depends only on the (boundary fixing) 
isotopy class of $\varphi$ and that 
this construction gives a monoid homomorphism 
from the mapping class group $\mathcal{M}_{g,1}$ 
of $\Sigma_{g,1}$ to $\mathcal{C}_{g,1}$. In fact, it is 
an injective homomorphism 
(see Garoufalidis-Levine \cite[Section 2.4]{gl}, Levine 
\cite[Section 2.1]{levin} 
and also \cite[Proposition 2.3]{gs}). 
It is known that the image of this homomorphism coincides with 
the group of units of $\mathcal{C}_{g,1}$ 
(see Habiro-Massuyeau \cite[Proposition 2.4]{hm} for example). 
\end{example}

As seen in Example \ref{ex:mgtocg}, 
we may regard $\Cg$ as an enlargement of $\Mg$. 
We recall here that $\Mg$ is a perfect group for $g \ge 3$ 
(see Korkmaz's survey \cite{korkmaz} and 
papers listed there for details). 
In comparing the structures of $\Cg$ and $\Mg$, 
it seems interesting to discuss abelian quotients of $\Cg$. 
Note that we need to be careful when mentioning 
abelian quotients of $\Cg$ since it 
is not a group but a monoid. 
We avoid this by considering its 
universal group $\mathcal{U}(\Cg)$ and its abelian quotients. 
Recall that for every monoid $S$, 
there exist (uniquely up to isomorphism) 
a group $\mathcal{U}(S)$ together with 
a monoid homomorphism $g:S \to \mathcal{U}(S)$ satisfying 
the following: For every monoid homomorphism 
$f: S \to G$ to a group $G$, 
there exists a unique group homomorphism 
$f':\mathcal{U}(S) \to G$ such that $f=f' \circ g$. 
One of the possible constructions 
of $\mathcal{U}(S)$ is to regard a monoid presentation of $S$ as a group 
presentation. 

In our discussion below, the case where $g=0$ is exceptional. 
We can check that $\mathcal{C}_{0,1}$ is isomorphic to 
the monoid $\theta_\Z^3$ of all (integral) homology 3-spheres 
whose product is given by connected sums. Indeed, 
an isomorphism $\theta_\Z^3 \cong \mathcal{C}_{0,1}$ 
is given by assigning to each homology 3-sphere $X$ the 
homology cylinder $((D^2 \times [0,1])\sharp X, \Id \times 1, \Id \times 0)$ 
with the inverse homomorphism given by taking {\it closures} 
(see Example \ref{ex:knot}). Consequently, 
$\mathcal{C}_{0,1}$ is an abelian monoid which is not finitely generated. 

We begin our main argument by the following observations. 
\begin{lemma}\label{lem:incompressible}
For each $(M,i_+,i_-) \in \Cg$ with $g \ge 1$, the surfaces 
$i_+(\Sg)$ and $i_-(\Sg)$ are incompressible in $M$.
\end{lemma}
\begin{proof}
It suffices to show that $i_+:\pi_1(\Sg) \to \pi_1(M)$ is injective. 
Since $i_+:\Sg \hookrightarrow M$ induces an isomorphism on homology, 
it follows from Stallings' theorem \cite[Theorem 3.4]{st} that $i_+$ induces 
isomorphisms on all stages of nilpotent quotients. Combining it with the 
fact that $\pi_1(\Sg)$ is free, in particular residually nilpotent, 
we see that $i_+:\pi_1(\Sg) \to \pi_1 (M)$ is injective. 
\end{proof}
\begin{proposition}\label{prop:notFG1}
The monoid $\mathcal{C}_{g,1}$ is {\rm not} finitely 
generated, for every $g \ge 0$. 
In fact, the abelianization of $\mathcal{U}(\Cg)$ has infinite rank.
\end{proposition}
\begin{proof}
The case where $g=0$ is as mentioned above. 
We now assume $g \ge 1$. 
For each homology cylinder $(M,i_+,i_-) \in \Cg$, 
the underlying 3-manifold $M$ has a prime decomposition of the form
\[M \cong M_0 \sharp X_1 \sharp X_2 \sharp \cdots \sharp X_n,\]
where $M_0$ is the unique prime factor containing $\partial M$ and 
$X_1, X_2, \ldots, X_n$ are homology 3-spheres. 
Using this decomposition, we can define the {\it forgetting} map 
\[F:\mathcal{C}_{g,1} \longrightarrow \theta_\Z^3\]
by
\[F(M,i_+,i_-) = S^3 \sharp X_1 \sharp X_2 \sharp 
\cdots \sharp X_n.\]
The uniqueness of the prime decomposition of a 3-manifold shows 
that $F$ is well-defined. 

We now claim that $F$ is a surjective monoid homomorphism. 
Let $(M,i_+,i_-), (N,j_+,j_-)$ $\in \Cg$. We decompose $M$ 
into $M \cong M_0 \sharp X_M$, where 
$M_0$ is the prime factor containing $\partial M$ 
and $X_M$ is a homology 3-sphere. 
Similarly, we have $N \cong N_0 \sharp X_N$. 
By Lemma \ref{lem:incompressible}, the underlying 3-manifold of 
the product $(M,i_+,i_-)\cdot(N,j_+,j_-)$ has a decomposition 
\[(M_0 \cup_{i_- \circ (j_+)^{-1}} N_0) \sharp X_M \sharp X_N\]
such that $M_0 \cup_{i_- \circ (j_+)^{-1}} N_0$ is the prime factor 
containing the boundary. This shows that $F$ is a monoid homomorphism. 
The surjectivity of $F$ follows from the existence of 
a section $\theta_\Z^3 \to \Cg$ 
defined by $X \mapsto ((\Sg \times [0,1]) \sharp X, 
\Id \times 1, \Id \times 0)$. 

The result follows from the fact that $\theta_\Z^3$ satisfies the conditions 
mentioned in the statement. 
\end{proof}
Proposition \ref{prop:notFG1} says that the monoid $\mathcal{C}_{g,1}$ 
has a different property about its abelian quotients 
from the mapping class group 
$\mathcal{M}_{g,1}$, which 
arises from the reducibility of the underlying 
3-manifold of a homology cylinder. However the underlying 3-manifolds 
of homology cylinders obtained from $\mathcal{M}_{g,1}$ are all 
product $\Sg \times [0,1]$ and, in particular, irreducible. 
Therefore it seems reasonable to consider the following subset of 
$\mathcal{C}_{g,1}$. 

\begin{definition}\label{def:irred}
A homology cylinder $(M,i_+,i_-)$ is said to be {\it irreducible} if 
the underlying 3-manifold $M$ is irreducible. We denote by 
$\mathcal{C}_{g,1}^{\mathrm{irr}}$ the subset of $\mathcal{C}_{g,1}$ 
consisting of all irreducible homology cylinders.
\end{definition}
\noindent
Note that $\mathcal{C}_{g,1}^{\mathrm{irr}}$ is a submonoid of 
$\mathcal{C}_{g,1}$, for 
$\mathcal{C}_{g,1}^{\mathrm{irr}} = \mathrm{Ker}(F)$. 
In particular, $\mathcal{C}_{0,1}^{\mathrm{irr}}$ is the trivial monoid. 
Note also that we have an injective monoid homomorphism 
$\mathcal{M}_{g,1} \hookrightarrow \mathcal{C}_{g,1}^{\mathrm{irr}}$. 
The following is the main result of this paper, whose proof 
will be given in the next section. 
\begin{theorem}\label{thm:notFG2}
The monoid $\mathcal{C}_{g,1}^{\mathrm{irr}}$ 
is not finitely generated, for every $g \ge 1$. 
In fact, the abelianization of $\mathcal{U}(\CgI)$ 
has infinite rank.
\end{theorem}

\section{Proof of Theorem \ref{thm:notFG2}}\label{sec:proof}
Our proof of Theorem \ref{thm:notFG2} will be obtained 
as an application of sutured Floer homology theory due to 
Juh\'asz \cite{juhasz,juhasz2} and Ni \cite{ni0,ni}.

For each homology cylinder $(M,i_+,i_-) \in \Cg$, 
we have a natural decomposition 
\[\partial M = i_+(\Sg) \cup_{i_+(\partial \Sg)=i_-(\partial \Sg)} 
i_-(\Sg)\]
of $\partial M$. Such a decomposition 
defines a {\it sutured manifold} $(M,\gamma)$ with the suture 
$s(\gamma) =i_+(\partial \Sg)=i_-(\partial \Sg)$. 
Sutured manifolds were originally 
defined by Gabai \cite{gabai1}, to which we refer for details. 
\begin{example}\label{ex:knot}
For a knot $K$ in a closed oriented connected 3-manifold $Y$ 
with a Seifert surface $S$, 
let $M$ be the manifold obtained from the knot exterior $\overline{Y - N(K)}$ 
by cutting open along $S$ and 
$\gamma$ the annulus $\partial N(K)\cap \partial M$.
Then, $(M,\gamma)$ is called the 
{\it complementary sutured manifold for} $S$.  
The core curve of $\gamma$ is denoted by 
$s(\gamma)$, and called the {\it suture\/}. 
The suture $s(\gamma)$ and $K$ are parallel. 
On the other hand, 
for each $(M,i_+,i_-) \in \Cg$ 
(or, more generally, each marked cobordism 
of $\Sg$), we have a closed oriented connected $3$-manifold 
\[C_M:=M/(i_+(x)=i_-(x)) \qquad (x \in \Sg)\]
called the {\it closure} and a knot $i_+(\partial \Sg)=i_-(\partial \Sg)$ 
with a Seifert surface $S=i_+(\Sg)=i_-(\Sg)$ in $C_M$. 
\end{example}

Sutured Floer homology 
is an invariant of {\it balanced} sutured manifolds 
(see Juh\'asz \cite[Definition 2.11]{juhasz}), 
where all of the sutured manifolds mentioned above 
satisfy this condition. 
It assigns a finitely generated abelian group $SFH(M,\gamma)$ 
to each balanced sutured manifold $(M,\gamma)$. We 
rely on papers of Juh\'asz \cite{juhasz,juhasz2} 
and Ni \cite{ni0,ni} for the definition and fundamental 
properties of sutured Floer homology, 
and concentrate on using this theory. 
Juh\'asz \cite[Theorem 1.5]{juhasz2} showed 
that 
\begin{equation}\label{eq:SFH-HFK}
SFH(M,\gamma) = \widehat{HFK}(C_M, s(\gamma), g(S))
\end{equation}
holds for any sutured manifold $(M,\gamma)$ mentioned in 
Example \ref{ex:knot}, where the right hand side is 
the genus $g(S)$ part of the knot Floer homology of 
the knot $s(\gamma)$  in $C_M$ 
with the Seifert surface $S$ of genus $g(S)$. 

For each homology cylinder $(M,i_+,i_-) \in \Cg$, 
we put $SFH(M,i_+,i_-):=SFH(M,\gamma)$ with 
$s(\gamma)=i_+(\partial \Sg)=i_-(\partial \Sg)$. 

\begin{proposition}\label{prop:nontrivial}
$SFH(M,i_+,i_-)$ contains \,$\Z$ \ for any $(M,i_+,i_-) \in \Cg$. 
\end{proposition}
\begin{proof}
We first assume that $M$ is irreducible. 
By Juh\'asz \cite[Theorem 1.4]{juhasz2}, 
all we have to do is 
to check that $(M,i_+,i_-)$ gives a {\it taut} sutured manifold. That is, 
\begin{itemize}
\item[(i)] $M$ is irreducible; and 
\item[(ii)] $i_+(\Sg)$ and $i_-(\Sg)$ are incompressible and Thurston 
norm minimizing in their homology classes in $H_2(M,\gamma)$ 
with $s(\gamma)=i_+(\partial \Sg)=i_-(\partial \Sg)$. 
\end{itemize}
The condition (i) is automatic and 
the first half of (ii) follows from Lemma \ref{lem:incompressible}. 
For the latter half of (ii), it suffices to 
show that $i_+ (\Sg)$ is Thurston norm minimizing. 

Suppose that we have a proper embedding of a surface $\Sigma$ 
representing $[i_+ (\Sg)] \in H_2 (M,\gamma)$ with 
smaller norm than that of $i_+ (\Sg)$. 
We may assume that $\Sigma$ does not 
have any closed component, for such a component comes from 
$H_2 (M)=0$ and hence removing it does not change 
the class $[\Sigma] \in H_2 (M,\gamma)$ nor increase the norm. 
Next, consider the intersection $\Sigma \cap \gamma$ 
of $\Sigma$ and the annulus $\gamma$. 
Generically, it consists of oriented circles, each of which 
is an essential simple loop in $\gamma$ or bounds a disk 
in $\gamma$. 
For an innermost circle bounding a disk, we attach this disk 
to $\Sigma$ and move it away from $\partial M$ by an isotopy. Repeating this, 
we can eliminate all circles bounding disks in $\Sigma \cap \gamma$, 
so that the intersection consists of parallel copies 
of the essential simple loop in $\gamma$. 
If $\Sigma \cap \gamma$ is disconnected, we can find an annulus in $\gamma$ 
bounding (with coherent orientations) 
two adjacent components of 
$\Sigma \cap \gamma$. Then we attach this annulus to 
$\Sigma$ and move it away from $\partial M$. 
(When a closed component is produced, remove it.) 
Repeating this, $\Sigma \cap \gamma$ becomes connected. 
Note that the above procedure does not change $[\Sigma] 
\in H_2 (M,\gamma)$ 
nor increase the norm since disks and annuli in $\gamma$ 
are trivial in $H_2 (M,\gamma) \cong H_1 (\gamma) \cong \Z$. 
Consequently, we may assume that $\Sigma$ is a connected surface 
with one boundary component. 

Now suppose that 
we have a proper embedding $j:\Sigma_{h,1} \hookrightarrow M$ 
with 
$h<g$ and satisfying $[j(\Sigma_{h,1})]=[i_+(\Sg)] 
\in H_2(M,\gamma)$. 
We take a basis 
$\{\delta_1, \delta_2, \ldots, \delta_{2h}\}$ of $\pi_1(j(\Sigma_{h,1}))$ 
such that 
\begin{equation}\label{eq:boundary}
s(\gamma)=i_+(\partial\Sg)=j(\partial\Sigma_{h,1})=
\prod_{i=1}^h [\delta_{2i-1}, \delta_{2i}] \in [\pi_1(M),\pi_1(M)]
\end{equation}
under suitable orientations. 
Here we set the basepoint on $s(\gamma)$.
By Stallings \cite[Lemma 3.1]{st}, 
we see that 
\[i_+:\frac{[\pi_1 (\Sg),\pi_1(\Sg)]}{[\pi_1 (\Sg),[\pi_1 (\Sg),\pi_1(\Sg)]]} 
\longrightarrow 
\frac{[\pi_1 (M),\pi_1(M)]}{[\pi_1 (M),[\pi_1 (M),\pi_1(M)]]}\]
is an isomorphism, and we pull back (\ref{eq:boundary}) to 
$\frac{[\pi_1 (\Sg),\pi_1(\Sg)]}{[\pi_1 (\Sg),[\pi_1 (\Sg),\pi_1(\Sg)]]}$, 
which is known to be isomorphic to $\wedge^2 (H_1(\Sg))$. Then we obtain 
an equality 
\[\partial\Sg=\sum_{i=1}^h \, i_+^{-1}([\delta_{2i-1}]) \wedge 
i_+^{-1}([\delta_{2i}]) 
\in \wedge^2 (H_1(\Sg)), \]
where $[\delta_k] \in H_1(M)$ denotes the homology class of $\delta_k$. 
On the other hand, we have 
$\partial \Sg=\sum_{j=1}^g x_j \wedge y_j$, the symplectic form, for 
any symplectic basis $\{x_1,\ldots,x_g,y_1,\ldots,y_g\}$ of $H_1(\Sg)$. 
Define an endomorphism of $H_1(\Sg)$ by 
\begin{align*}
&x_i \mapsto i_+^{-1}([\delta_{2i-1}]), \ \ 
y_i \mapsto i_+^{-1}([\delta_{2i}]) \quad \mbox{for \ $1 \le i \le h$},\\
&x_j, y_j \mapsto 0 \quad \mbox{for \ $h+1 \le j \le g$}. 
\end{align*}
By definition, this endomorphism is not injective, but 
preserves the symplectic form. However, 
such an endomorphism does not exist 
since the symplectic form embodies the intersection 
form on $H_1 (\Sg)$, which is nondegenerate, a contradiction. 
Therefore $i_+ (\Sg)$ is Thurston norm minimizing and we finish the proof 
when $M$ is irreducible. 

When $M$ is not irreducible, we take a prime decomposition 
$M=M_0 \sharp X_1 \sharp X_2 \sharp \cdots \sharp X_n$ 
as in Section \ref{section:cylinder}, where 
$X_1, X_2, \ldots, X_n$ are all homology 3-spheres. Then 
we obtain the conclusion by an argument similar to the proof of 
\cite[Corollary 8.3]{juhasz2} using the connected sum 
formula \cite[Proposition 9.15]{juhasz}. 
\end{proof}

By formulas of Juh\'asz 
\cite[Proposition 8.6]{juhasz2} and 
Ni \cite[Theorem 4.1, 4.5]{ni} together 
with the fact (\ref{eq:SFH-HFK}), 
we have 
\[SFH((M,i_+,i_-) \cdot (N,j_+,j_-))\otimes \Q \cong 
(SFH(M,i_+,i_-) \otimes SFH(N,j_+,j_-)) \otimes \Q\]
for $(M,i_+,i_-)$, $(N,j_+,j_-) \in \Cg$. 
Hence by taking the rank of $SFH$, 
we obtain a monoid homomorphism
\[R: \Cg \longrightarrow \Z_{>0}^\times\]
defined by 
\[R(M,i_+,i_-) =  \mathrm{rank}_\Z (SFH(M,\gamma)),\]
where 
$\Z_{>0}^\times$ is the monoid 
of positive integers 
whose product is given by multiplication. 
We call $R$ the {\it rank homomorphism}. 
Note that the restriction of $R$ to $\Mg$ is trivial 
since every element of $\Mg$ has its inverse. 

By the uniqueness of the prime decomposition of an integer, 
we can decompose $R$ into prime factors
\[R=\bigoplus_{\scriptsize \mbox{$p$\,: prime}} 
R_p : \Cg \longrightarrow 
\Z_{>0}^\times = 
\bigoplus_{\scriptsize \mbox{$p$\,: prime}} \PZ^{(p)},\]
where $\PZ^{(p)}$ is a copy of $\PZ$, 
the monoid of non-negative integers 
whose product is given by sums, 
corresponding to the power of the prime number $p$. 
We now restrict the above homomorphisms to $\CgI$. 
\begin{proposition}\label{prop:independent}
For $g \ge 1$, the set 
$\{R_p:\CgI \to \PZ \mid p\ \mathrm{prime}\}$
contains infinitely many non-trivial homomorphisms 
that are linearly independent as elements of $\Hom (\CgI,\PZ)$. 
\end{proposition}
To prove this proposition, 
we have to consider the images of the homomorphisms $R_p$. 
More specifically, we need many 
homology cylinders whose ranks of $SFH$ are known. We now 
use {\it homologically fibered knots} 
defined in the previous paper \cite{gs} to construct 
such homology cylinders. 

Recall that a knot $K$ in $S^3$ is said to be 
{\it homologically fibered} if 
$K$ satisfies the following two conditions:
\begin{itemize}
\item[(i)] The degree of the normalized 
Alexander polynomial $\Delta_{K}(t)$ of $K$ 
is $2g(K)$; and 
\item[(ii)] $\Delta_{K}(0)=\pm 1$,
\end{itemize}
\noindent
where $g(K)$ is the genus of $K$ and 
$\Delta_{K}(t)$ is normalized 
so that its lowest degree is $0$. 
In \cite[Theorem 3.4]{gs}, we showed that 
$K$ is homologically fibered if and only if 
$K$ has a Seifert surface $S$ 
whose complementary sutured manifold is 
a homology product. Here, a homology product means 
a homology cylinder without markings. 
Note that Crowell and Trotter
observed in \cite{ct} this essentially. 
(See also \cite{ni}.) 
We also showed that if $K$ is homologically fibered, then 
any minimal genus Seifert surface gives a homology product. 

For a homologically fibered knot $K$ of genus $g$, 
we obtain an irreducible homology cylinder 
$(M,i_+,i_-) \in \CgI$ by fixing a pair of markings of 
the boundary of the complementary sutured manifold $(M,K)$ 
for a minimal genus Seifert surface. 
By (\ref{eq:SFH-HFK}), 
\[SFH(M,i_+,i_-)=SFH(M,K) \cong \widehat{HFK} (S^3,K,g)\]
holds for such a homology cylinder. 

\begin{proof}[Proof of Proposition $\ref{prop:independent}$]
We first give a proof of the case where $g=1$. 
We now consider 
pretzel knots $P(2l+1,2m+1,2n+1)$ with $2l+1<0$. As depicted in 
Figure \ref{pretzel1}, each of such knots has a 
genus $1$ Seifert surface. 

\begin{figure}[h]
\centering
\includegraphics[width=.8\textwidth]{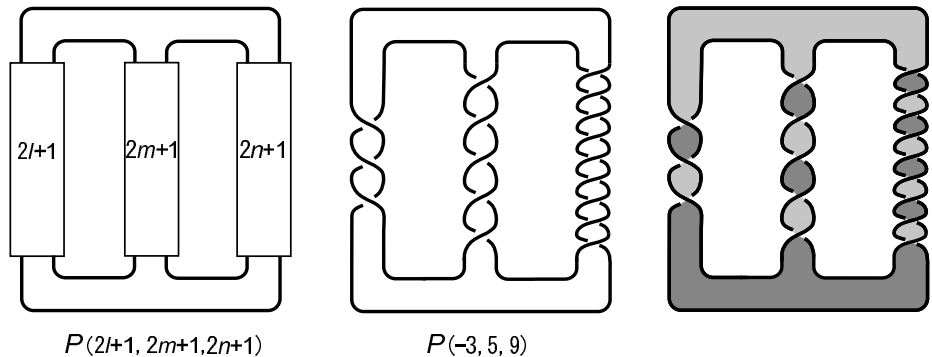}
\caption{Standard diagram of Pretzel knots with Seifert surfaces 
of genus 1}\label{pretzel1}
\end{figure}

It is well known that 
the normalized Alexander polynomial of $P(2l+1,2m+1,2n+1)$ 
is given by 
\[\Delta_{P(2l+1,2m+1,2n+1)}(t)=
(1+l+m+n+lm+mn+nl)(t-1)^2+t.\]
Using this formula, we see that the sequence 
\[\{P_n:=P(-2n+1, 2n+1, 2n^2+1)\}_{n=1}^\infty\]
consists of homologically fibered knots of genus $1$ since 
$\Delta_{P_n}(t)=1 -t +t^2$. Moreover, a computation due to 
Ozsv\'ath-Szab\'o \cite[Latter formula of Theorem 1.3]{os} 
gives
\[\widehat{HFK}(S^3,P_n,1) \cong \Z_{(1)}^{n^2-n} \oplus 
\Z_{(2)}^{n^2-n+1} \cong \Z^{2n^2-2n+1},\]
from which we see that $\{2n^2-2n+1\}_{n=1}^\infty 
\subset R(\mathcal{C}_{1,1}^{\mathrm{irr}})$. 

We now analyze the sequence $\{2n^2-2n+1\}_{n=1}^\infty$ 
of positive integers. 
The first supplement to quadratic reciprocity 
says that there exists 
a positive integer $m$ such that $m^2 \equiv -1 \bmod p$ 
for any odd prime $p$ satisfying $p \equiv 1 \bmod 4$. 
In this case, set 
\[n=\left\{\begin{array}{cl}
\displaystyle\frac{m+1}{2} & \mbox{($m$: odd)} \\ &\\
\displaystyle\frac{m+p+1}{2} & \mbox{($m$: even)}
\end{array}\right. .\]
Then, if $m$ is odd, we have 
\[2n^2-2n+1=\frac{m^2+1}{2}\equiv 0 \quad \bmod p,\]
and if $m$ is even, we also have 
\[2n^2-2n+1=\frac{(m+p)^2+1}{2} \equiv 0 \quad \bmod p.\]
Hence we can conclude that $R_p$ is non-trivial if 
$p \equiv 1 \bmod 4$. 
By Dirichlet's theorem on arithmetic progressions, 
there exist 
infinitely many such prime numbers. 

For a homology cylinder $M \in \mathcal{C}_{1,1}^{\mathrm{irr}}$, 
let 
\[p(M) := \max(\{1\} \cup \{p \mid R_p (M)\neq 0\}).\]
Take a sequence of homology cylinders $\{M_i\}_{i=1}^\infty 
\subset \mathcal{C}_{1,1}^{\mathrm{irr}}$ such that 
$p(M_1) < p(M_2) < \cdots$. Then we can see that 
$\{R_{p(M_i)}\}_{i=1}^\infty$ are linearly independent by evaluating 
them on the $M_i$'s. 
This concludes the proof of the case where $g=1$. 

Let $P_n(k)$ be a homologically fibered knot of 
genus $k+1$ obtained from 
$P_n$ by taking connected sums with $k$-tuples of trefoils. 
By \cite[Theorem 1.1]{ni0} and \cite[Corollary 8.8]{juhasz2} together with 
the fact that the trefoil is a fibered knot of genus $1$, so that 
$\widehat{HFK}(S^3,\mathrm{trefoil},1) \cong \Z$, we have 
\[\mathrm{rank}_\Z (\widehat{HFK}(S^3,P_n,1)) = 
\mathrm{rank}_\Z (\widehat{HFK}(S^3,P_n(k),k+1)).\]
Therefore, the cases where $g \ge 2$ follow from 
the same argument as above. 
\end{proof}

\begin{proof}[Proof of Theorem $\ref{thm:notFG2}$]
Suppose $\CgI$ was finitely generated. Then except for 
finitely many primes, the homomorphisms $R_p$ are trivial 
on any finite set of generators, and hence on whole $\CgI$. 
This contradicts Proposition \ref{prop:independent} and we have proved 
the first half of our claim. The latter half follows from 
the construction that uses infinitely many homomorphisms whose targets 
are abelian. 
\end{proof}

\section{Observations from the viewpoint of homology cobordism}
\label{sec:Hcob}
In \cite{gl}, Garoufalidis-Levine introduced {\it homology cobordisms} 
of homology cylinders, 
which give an equivalence relation of homology cylinders. 
We finish this paper by two observations concerning our results 
(Proposition \ref{prop:notFG1} and Theorem \ref{thm:notFG2}) from 
the viewpoint of this equivalence relation. 

Two homology cylinders $(M,i_+,i_-)$ and $(N,i_+,i_-)$ over 
$\Sigma_{g,1}$ are {\it homology cobordant} 
if there exists a compact oriented smooth 4-manifold $W$ such that: 
\begin{enumerate}
\item $\partial W = M \cup (-N) /(i_+ (x)= j_+(x) , \,
i_- (x)=j_-(x)) \quad x \in \Sigma_{g,1}$; 
\item the inclusions $M \hookrightarrow W$, $N \hookrightarrow W$ 
induce isomorphisms on the integral homology.  
\end{enumerate}
\noindent
We denote by $\mathcal{H}_{g,1}$ 
the quotient set of $\mathcal{C}_{g,1}$ with respect 
to the equivalence relation of homology cobordism. 
The monoid structure of $\mathcal{C}_{g,1}$ induces 
a group structure of $\mathcal{H}_{g,1}$. 
It is known that $\Mg$ can be embedded in $\Hg$ 
(see \cite[Section 2.4]{gl}, \cite[Section 2.1]{levin}). 

One important problem is to determine whether $\Hg$ is perfect or not
\footnote{After the authors wrote the first version of this paper, 
Cha, Friedl and Kim settled this problem in \cite{cfk}, where 
they showed that $\Hg$ has $(\Z/2\Z)^\infty$ as an abelian quotient, 
for every $g \ge 1$.}. 
In fact, no non-trivial abelian quotients of $\mathcal{H}_{g,1}$ 
are known at present. We now observe that it is difficult 
to give an answer to this problem by using the homomorphisms used 
in this paper. First we consider the forgetting homomorphism 
$F: \Cg \to \theta_\Z^3$ discussed 
in Section \ref{section:cylinder}. 

\begin{theorem}\label{thm:ab_quot1}
For every abelian group $A$ and every non-trivial monoid homomorphism 
$\varphi_A:\theta_\Z^3 \to A$, 
the composite $\varphi_A \circ F:\mathcal{C}_{g,1} \to A$ 
does not factor 
through $\mathcal{H}_{g,1}$, for all $g \ge 1$. 
\end{theorem}
\begin{proof}
It follows from 
Myers' result \cite[Theorem 3.2]{myers} that 
every homology cylinder in $\Cg$ with $g \ge 1$ 
is homology cobordant to 
an irreducible one, whose image by $F$ is 
trivial by definition. 
Hence if $\varphi_A \circ F$ factors through $\Hg$ for a 
monoid homomorphism $\varphi_A: \theta_\Z^3 \to A$, then 
$\varphi_A$ must be trivial. 
\end{proof}

Next we consider the rank homomorphisms $R_p : \CgI \to \PZ$ 
discussed in Section \ref{sec:proof}. It induces 
a group homomorphism 
$R_p:\mathcal{U}(\mathcal{C}_{g,1}^{\mathrm{irr}}) \to 
\Z$ on universal groups. 
Note that the quotient group of $\Cg^{\mathrm{irr}}$ by 
homology cobordism relation is also $\Hg$ as mentioned 
in the proof of Theorem \ref{thm:ab_quot1}. 
\begin{theorem}\label{thm:ab_quot2}
For each $g \ge 1$, the homomorphism 
$R_p : \mathcal{U}(\mathcal{C}_{g,1}^{\mathrm{irr}}) \to 
\Z$ does not factor through $\mathcal{H}_{g,1}$ if it is non-trivial. 
\end{theorem}
\begin{proof}
Since $R_p(\CgI) \subset \PZ$ and 
$\Hg$ is a quotient group of $\CgI$, 
the homomorphism $R_p$ must be trivial if it factors through $\Hg$. 
\end{proof}

\bibliographystyle{amsplain}

\end{document}